\documentclass[11 pt]{amsart}
\usepackage{amsfonts}
\usepackage{ifthen}
\usepackage{amsthm}
\usepackage{amsmath}
\usepackage{graphicx}
\usepackage{amscd,amssymb,amsthm}
\usepackage{graphicx}
\usepackage{epstopdf}
\usepackage{hyperref}
\usepackage{cleveref}
\usepackage{cite}
\usepackage{mathrsfs}

\newcounter{minutes}
\divide\time by 60
\newcounter{hours}
\multiply\time by 60 \addtocounter{minutes}{-\time}


\begin{document}
	
	\title{On some geometric properties of normalized Wright functions}

	\author[E. Toklu]{Evr{\.I}m Toklu}
	\address{Department of Mathematics, Faculty of Education, A\u{g}r{\i} {\.I}brahim \c{C}e\c{c}en University, A\u{g}r{\i}, Turkey} 
	\email{evrimtoklu@gmail.com}
	
	\author[N. Karag\"{o}z]{Nesl{\.I}han Karag\"{O}z}
	\address{Institute of Science and Technology, A\u{g}r{\i} {\.I}brahim \c{C}e\c{c}en University, A\u{g}r{\i}, Turkey} 
	\email{nslhnkrgz@gmail.com}
	
	\def\thefootnote{}
	\footnotetext{ \texttt{File:~\jobname .tex,
			printed: \number\year-0\number\month-\number\day,
			\thehours.\ifnum\theminutes<10{0}\fi\theminutes}
	} \makeatletter\def\thefootnote{\@arabic\c@footnote}\makeatother
	
	\maketitle
	
\begin{abstract}
The main purpose of the present paper is to determine the radii of starlikeness and convexity associated with lemniscate of Bernoulli and the Janowski function, $(1+Az)/(1+Bz)$ for $-1\leq B<A\leq 1,$ of normalized Wright functions. The key tools in the proof of our main results are the infinite product representation of Wright function and properties of real zeros of the Wright function and its derivative.
\end{abstract}
\section{\bf Introduction and the main results}
Let $\mathbb{D}_{r}$ be the open disk $\left\lbrace z\in\mathbb{C}:\left| z\right| <r\right\rbrace $ with the radius $r>0$ and let $\mathbb{D}=\mathbb{D}_{1}.$ Let $f:\mathbb{D}_r\to\mathbb{C}$ be the function defined by
\begin{equation}
f(z)=z+\sum_{n\geq 2}a_{n}z^{n},  \label{W0}
\end{equation}
where $r$ is less or equal than the radius of convergence of the above power series. Let $\mathcal{A}$ be the class of analytic functions of the form \eqref{W0}, that is, normalized by the conditions $f(0)=f^{\prime}(0)-1=0.$ Let $\mathcal{S}$ denote the class of functions belonging to $\mathcal{A}$ which are of univalent in $\mathbb{D}_{r}$. A function $f\in \mathcal{A}$ is said to be starlike function if $f(\mathbb{D})$ is starlike domain with the respect to the origin. It is well known fact that various subclasses of starlike function can be unified by making use of the concept of subordination. A function $f\in \mathcal{A}$ is subordinate to a function $g\in \mathcal{A},$ written as $f(z)\prec g(z),$ if there exist a Schwarz function $w$ with $w(0)=0$ and $\left| w(z)\right| <1$ such that $f(z)=g(w(z)).$ In addition, we know that if $g$ is a univalent function, then $f(z)\prec g(z)$ if and only if $f(0)=g(0)$ and $f(\mathbb{D})\subset g(\mathbb{D}).$ For an analytic function $\varphi,$ let $\mathcal{S}^{\star}(\varphi)$ denote the class of all analytic functions satisfying $1+zf'(z)/f(z)\prec \varphi(z).$ By $\mathcal{K}(\varphi)$ we mean the class of all analytic functions satisfying $1+zf''(z)/f'(z)\prec \varphi(z).$ It is worth mentioning that these classes include respectively several famous subclasses of starlike and convex functions. For instance, the class $\mathcal{S}^{\star}_{\mathcal{L}}:=\mathcal{S}^{\star}(\sqrt{1+z})$ denotes the class of lemniscate starlike functions introduced and investigated by Sok\'ol and Stankiewich \cite{SS} and the class $\mathcal{K}_{\mathcal{L}}:=\mathcal{K}(\sqrt{1+z})$ represents the class  of lemniscate convex functions. Moreover, for $-1\leq B<A\leq 1,$ the class $\mathcal{S}^{\star}[A,B]:=\mathcal{S}^{\star}((1+Az)/(1+Bz))$ is the class of Janowski starlike functions and $\mathcal{K}[A,B]:=\mathcal{K}((1+Az)/(1+Bz))$ is the class of Janowski convex functions \cite{Janowski}.

Given a class of functions $\mathcal{M}\subset\mathcal{A}$ and a function $f\in \mathcal{A},$ the $\mathcal{M}-$radius of the function $f$ is the largest number $r$ with $0\leq r \leq 1$ such that $f_{r}\in \mathcal{M},$ where $f_{r}(z):=f(rz)/r.$ If we choose $\mathcal{M}=\mathcal{S}^{\star}_{\mathcal{L}},$ the $\mathcal{M}-$radius of the function $f$, which is represented by $r^{\star}_{\mathcal{L}}(f),$ is called the radius of lemniscate starlikeness. It is indeed the largest $r$ with $0\leq r \leq 1$ such that
\[\left|\left(\frac{zf'(z)}{f(z)} \right)^2-1  \right|<1 \quad (\left| z\right|<r ).\]
If we choose $\mathcal{M}=\mathcal{K}_{\mathcal{L}},$ the $\mathcal{M}-$radius of the function $f$, which is represented by $r^{c}_{\mathcal{L}}(f),$ is called as the radius of lemniscate convexity. It is indeed the largest $r$ with $0\leq r \leq 1$ such that
\[\left|\left(1+\frac{zf''(z)}{f'(z)} \right)^2-1  \right|<1 \quad (\left| z\right|<r ).\]
If we take $\mathcal{M}=\mathcal{S}^{\star}[A,B]$ or $\mathcal{M}=\mathcal{K}[A,B],$ the respective $\mathcal{M}-$radii, which are represented by $r^{\star}_{A,B}(f)$ and $r^{c}_{A,B}(f),$ are called as the radii of Janowski starlikeness and Janowski convexity. These are respectively the largest  $r$ with $0\leq r \leq 1$ such that
\[ \left|\frac{\left( zf'(z)/f(z)\right)-1 }{A-Bzf'(z)/f(z)} \right|<1 \text{ \ and \ }      \left|\frac{zf''(z)/f'(z) }{A-B\left( 1+zf''(z)/f'(z)\right) }\right| <1 \quad (\left| z\right|<r ).\]

Recently, there has been a vivid interest on some geometric properties such as univalency,
starlikeness, convexity and uniform convexity of various special functions such as hyper-Bessel, Wright, $q-$Bessel and Mittag-Leffler functions (see \cite{ABskecul,ABS,ATO,AP,BTK,TAO}). Fore more details on the radius problems, one may consult on the studies \cite{AJR,Goodman,MKR1,Toklu,VR}. Moreover, in \cite{MKR2} the authors determined the radii of starlikeness and convexity associated with lemniscate of Bernoulli and the Janowski function $(1+Az)/(1+Bz)$.  Motivated by the above series of papers on geometric properties of special functions, in this paper our aim is to determine the radii of lemniscate starlikeness, lemniscate convexity, Janowski starlikeness and Janowski convexity of normalized Wright functions.

\end{document}